\documentclass[11pt]{amsart}
\usepackage{geometry}                
\usepackage{graphicx}
\usepackage{amssymb}
\usepackage{epstopdf}
\title{FRATTINI PROPERTIES AND NILPOTENCY IN LEIBNIZ ALGEBRAS}
\author{Allison McAlister, Kristen Stagg Rovira and Ernie Stitzinger}

\begin{document}
\maketitle

ABSTRACT

\vspace{5mm}

Ideals that share properties with the Frattini ideal of a Leibniz algebra are studied. 
Similar investigations have been considered in group theory. 
However most of the results are new for Lie algebras.
Many of the results involve nilpotency of these algebras.

\vspace{5mm}

I. INTRODUCTION

\vspace{5mm}

Frattini theory for algebras goes back at least 50 years. 
A general theory can be found in [14] and there are many works on special classes of algebras, especially Lie algebras. 
Leibniz algebras, as a generalization of Lie algebras, is a natural class to investigate and [2-8] contain results on Frattini subalgebras and ideals. 
Frattini theory for groups goes back to the nineteenth century, and there have been many results that are similar in groups and Lie algebras. 
Subgroups that have Frattini-like properties have been considered in [9] and related special types of subgroups have been studied in [11]. 
It is the purpose of this paper to view the Leibniz algebra analogues to some of these theories. 
Many of these results are new for Lie algebras as well. 
We consider only finite dimensional Leibniz algebras over a field F. 
The intersection of all maximal subalgebras of A is called the Frattini subalgebra of A and is denoted by F(A). 
Even in the solvable case, it need not be an ideal in A [5]. 
The maximal ideal of A contained in F(A) is called the Frattini ideal of A and is denoted by $\Phi$(A). 
References for Leibniz algebras include [1],[2],[10] and [12].

\vspace{5 mm}

II. GENERALIZED FRATTINI IDEALS

\vspace{5mm}

In [9], a proper subgroup, H, of a finite group G is called generalized Frattini if whenever G=HN$_G$(P) for any Sylow subgroup P of any normal subgroup K of G, then G=N$_G$(P). 
To consider such a property in Leibniz algebras, we replace Sylow subgroups with Cartan subalgebras. 
Unlike the group theory case, Frattini subalgebras do not have to be invariant, and we will use the Frattini ideal as our model to be formalized. 
To guarantee existence of Cartan subalgebras,  we assume the algebras are over an infinite field, [2]. 
In this context, an ideal H of A is generalized Frattini in A if whenever A=H+N$_A$(C), where C is a Cartan subalgebra of ideal K in A, it follows that A=N$_A$(C).  
We will show H is a generalized Frattini ideal of A if and only if whenever D and B are ideals of A and D is contained in B$\cap$H, then B/D nilpotent implies that B is nilpotent, a property possessed by Frattini ideals. 
We will find examples of this concept and conditions that guarantee that an ideal is generalized Frattini.

\vspace{5 mm}

Proposition 1. Let H be a generalized Frattini ideal in A. Then the following are true.

1. H is nilpotent.

2. Any ideal of A that is contained in H is also a generalized Frattini ideal in A

3. H+$\Phi$(A) is a generalized Frattini ideal in A

4. H+Z(A) is a generalized Frattini ideal in A whenever H+Z(A) is a proper subalgebra of A

\vspace{5 mm}

Proof. 1. Let C be a Cartan subalgebra of H. Then A=H+N$_A$(C) by Theorem 6.6 in [2]. H is generalized Frattini in A, hence A=N$_A$(C). Therefore, H=N$_H$(C). Since C is a Cartan subalgebra of H, N$_H$(C)=C. Thus H=C and H is nilpotent.

\vspace{5 mm}

2. Let N be an ideal of A such that N$\subseteq$ H. Let K be an ideal of A and let C be a Cartan subalgebra of K such that A=N+N$_A$(C). Then A=H+N$_A$(C) and hence, A=N$_A$(C) since H is generalized Frattini in A. Thus by definition, N is also generalized Frattini in A.

\vspace{5 mm}

3. Let K be an ideal in A with Cartan subalgebra C such that A=H+$\Phi$(A)+N$_A$(C). Suppose that M is a maximal subalgebra of A such that H+N$_A$(C) $\subseteq$ M. Then H+$\Phi$(A)+N$_A$(C) $\subseteq$ M, a contradiction. Hence A=H+N$_A$(C). Therefore A=N$_A$(C) since H is generalized Frattini in A. Therefore H+$\Phi$(A) is generalized Frattini in A.

\vspace{5mm}

4. Suppose that K is an ideal in A with Cartan subalgebra C such that A = H+Z(A)+N$_A$(C). Then A=H+N$_A$(C) and A=N$_A$(C) since H is generalized Frattini in A. Therefore H+Z(A) is also generalized Frattini in A.

\vspace{5mm}

Corollary 2.In a Leibniz algebra, A, both Z(A) and $\Phi$(A) are generalized Frattini in A.

\vspace{5mm}

Lemma 3. Any proper ideal, H, of a nilpotent Leibniz algebra A is generalized Frattini in A.

\vspace{5mm}

Proof. Let K be an ideal in A with Cartan subalgebra C such that H+N$_A$(C)=A. Then C=K and N$_A$(C)=A.

\vspace{5mm}

The next result shows that an important property of the Frattini ideal is shared with any generalized Frattini ideal.

\vspace{5mm}

Theorem 4. Let H be generalized Frattini in A. If K is an ideal in A that contains H and K/H is nilpotent, then K is nilpotent.

\vspace{5mm}

Proof. Let K be as in the statement of the theorem and let C be a Cartan subalgebra of K. Then (C+H)/H is a Cartan subalgebra of K/H Since K/H is nilpotent, K/H=(C+H)/H and K=C+H. Furthermore A=K+N$_A$(C) by Theorem 6.6 of [2]. Then A=K+N$_A$(C)=H+C+N$_A$(C)=H+N$_A$(C)=N$_A$(C) since H is generalized Frattini. Hence K=N$_K$(C)=C and K is nilpotent.

\vspace{5mm}
Corollary 5. Suppose that A is not 0. A is nilpotent if and only if A$^2$ is generalized Frattini. 

\vspace{5mm}

Proof. If A$^2$ is generalized Frattini, then the result follows from Theorem 4. If A is nilpotent, then the result follows from  Lemma 3.

\vspace{5mm}

Corollary 6. Let H be generalized Frattini in A. If K is an ideal in A  such that K$^{\omega} \subseteq$ H. Then K is nilpotent.

\vspace{5mm}

Proof. Let $\sigma$ be the natual mapping from K/K$^{\omega}$ onto K+H/H. K+H/H is nilpotent and then K+H is nilpotent by  Theorem 4.

\vspace{5mm}

Theorem 7 Let H be an ideal in a Leibniz algebra A. H is generalized Frattini in A if and only if for each ideal J of A that contains H, whenever J/H is nilpotent, then J is nilpotent.

\vspace{5mm}

Proof. If H is generalized Frattini, then the result is Theorem 4. Conversely, suppose that the condition on ideals J holds. Let K be an ideal of A, C a Cartan subalgebra of K with A=H+N$_A$(C). Then (C+H)/H is an ideal in A/H, hence also in (K+H)/H. Since (C+H)/H is Cartan in (K+H)/H, (C+H)/H=(K+H)/H. Therefore (K+H)/H is nilpotent and K+H is nilpotent by hypothesis. Therefore, K=C. Hence N$_A$(C)=N$_A$(K)=A and H is generalized Frattini in A.

\vspace{5mm}

Example 8. Let A be a Leibniz algebra with basis x,y,z and multiplications xz=x=-zx, zy=y=-yz, and xy=yx= 0. Let H=(x) and K=(y). H and K are generalized Frattini in A but H+K is not. Thus the sum of two generalized Frattini ideals need not be generalized Frattini. Note that this example is Lie so the result stands in Lie algebras as well.

\vspace{5mm}

Theorem 9. Let H be a generalized Frattini ideal in A and let K be an ideal of A that contains H. Then K/H is generalized Frattini in A/H if and only if K is generalized Frattini in A.

\vspace{5mm}

Proof. Suppose that K is generalized Frattini in A. Let J/H be an ideal in A/H such J/H contains K/H and (J/H)/(K/H) is nilpotent. Then J/K is nilpotent and, since K is generalized Frattini in A,  J is nilpotent by Theorem 7. Hence J/H is nilpotent and K/H is generalized Frattini in A/H by Theorem 7.

\vspace{5mm}

Conversely, suppose that K/H is generalized Frattini in A/H. Let J be an ideal in A which contains K such that J/K is nilpotent. Then (J/H)/(K/H) is nilpotent. Hence J/H is nilpotent since K/H is generalized Frattini, and then J is nilpotent since H is generalized Frattini. Therefore K is generalized Frattini in A by Theorem 7. 

\vspace{5mm}

Proposition 10. If Nil(A) is generalized Frattini in A, then every solvable ideal of A is nilpotent and is generalized Frattini in A.

\vspace{5mm}

Proof. Suppose that Nil(A) is generalized Frattini in A, let H be a solvable ideal in A and k be the smallest positive integer such that H$^{(k+1)}$=0. Then H$^{(k)}$ is abelian and H$^{(k)} \subseteq$ Nil(A).Then H$^{(k)}$ is generalized Frattini in A by Theorem 1. Working up the derived series in this manner, we find that H $\subseteq$ Nil(A). Thus H is nilpotent and it is contained in a generalized Frattini ideal, Nil(A). Hence H is generalized Frattini in A.

\vspace{5mm}

Corollary 11. If Nil(A) is generalized Frattini in A, then A is not solvable. 

\vspace{5mm}

Proof. A is not nilpotent since a generalized Frattini ideal is a proper ideal. If A is solvable, then it is nilpotent by Proposition 10, a contradiction.

\vspace{5mm}

Example 12. Continuing Example 8, Nil(A)=H+K is again seen to be not generalized Frattini since A is solvable.

\vspace{5mm}

Proposition 13. If H is generalized Frattini in A, then Nil(A/H)=Nil(A)/H.

\vspace{5mm}

Proof. Since H is nilpotent, H $\subseteq$ Nil(A). Nil(A)/H is clearly contained in Nil(A/H). Suppose that B is an ideal of A such that B/H=Nil(A/H). Then B is nilpotent by Theorem 4. Hence B is contained in Nil(A) and B/H is contained in Nil(A)/H.

\vspace{5mm}

Corollary 14. Let A be a non-nilpotent Leibniz algebra. Then Nil(A) is generalized Frattini in A if and only if Nil(A)=Rad(A), where Rad(A) is the maximal solvable ideal of A.

\vspace{5mm}

Proof: Suppose that Nil(A)=Rad(A). Let N be an ideal in A containing Nil(A) such that N/Nil(A) is nilpotent. Then N is solvable and Nil(A) $\subseteq$ N $\subseteq$ Rad(A)=Nil(A). Hence N is nilpotent and Nil(A) is generalized Frattini in A by Theorem 7. Conversely, suppose that Nil(A) is generalized Frattini in A. Then Rad(A) is nilpotent by Proposition 10 and Nil(A)=Rad(A).

\vspace{5mm}

Example 15. Let A=gl(n,F). Then Nil(A)=Rad(A)=Z(A) is generalized Frattini in A. 

\vspace{5mm}

In [2] and [3] Barnes extends his theory of Engel subalgebras from Lie to Leibniz algebras. For a$\in$ A, set E$_A$(a) be the Fitting null component of left multiplication by a on A. This space is a subalgebra called the Engel subalgebra for a. He notes that although a may not be in  E$_A$(a), there is a b in E$_A$(a) such that E$_A$(a)=E$_A$(b). Hence when working with these subalgebras, we usually can assume that a is in E$_A$(a). For a subalgebra U of A,if the Engel subalgebra for u in U both contains U and  is minimal in the set of Engel subalgebras for all elements in U,  then the Engel subalgebras for all elements in U contain U. He then shows C is a Cartan subalgebra of A if and only if C is minimal in the set of Engel subalgebras of A.

\vspace{5mm}

Theorem 16. Let H be an ideal in A. Then H is generalized Frattini in A if and only if for each ideal K of A and each Cartan subalgebra C of K, A=E$_A$(c) whenever A=H+E$_A$(c), for all c $\in$ C .

\vspace{5mm}

Proof. Let H be generalized Frattini in A. Let K, C be as in the theorem such that for each c$\in$C, A=H+E$_A$(c). (C+H)/H is a Cartan subalgebra of (K+H)/H by Theorem 6.3 of [2].  Then, using Engel's theorem, C acts nilpotently on A/H since A=H+E$_A$(c) for all c, hence also on (K+H)/H. Then (C+H)/H also acts nilpotently on (K+H)/H. Then (C+H)/H =(K+H)/H. Therefore (K+H)/H is nilpotent, as is K+H since H is generalized Frattini in A. Hence K is nilpotent and K=C. Since c$\in$ C=K and K is an ideal, A=E$_A$(c).

\vspace{5mm}

Conversely, suppose that H satisfies the conditions in the theorem. Let K be an ideal with H $\subseteq$ K with K/H nilpotent. Let C be a Cartan subalgebra of K and let c $\in$ C be an element such that E$_K$(c) is minimal in the set of Engel subalgebras for c $\in$C. Then A=K+E$_A$(c)= H+C+E$_A$(c) since K/H is nilpotent. Since C $\subseteq$ E$_A$(c), A=H+E$_A$(c). Hence A=E$_A$(c) by supposition and K= E$_K$(c)= C. Therefore K is nilpotent and H is generalized Frattini by Theorem 7.

\vspace{5mm}

\vspace{5mm}

III. PRIMITIVE IDEALS

\vspace{5mm}

An ideal K of A is primitive if 

1. $\Phi$(A/K)=0

2. A/K contains a unique minimal ideal

3. dim(A/K) $>$ 1

\vspace{5mm}

Example 17. Let A be the three dimensional cyclic Leibniz algebra generated by a with aa$^3$=a$^2$. Let K be the ideal with basis a$^2$+a$^3$. Then A/K has basis a,a$^2$ where we delete K from the notation. A/K is cyclic with generator a and aa$^2$=-a$^2$. The minimum polynomial for L$_a$ is x(x+1). Thus A/K has 2 maximal subalgeras and $\Phi$(A/K)=0 using section 4 of [5] and has a  unique minimal ideal. Hence K is a primitive ideal in A.

\vspace{5mm}

Lemma 18. Let K be a primitive ideal in A. Then K contains $\Phi$(A) and A/K is non nilpotent. Hence A is non-nilpotent.

\vspace{5mm}

Proof. Since $\Phi$(A/K)=0, $\Phi$(A) $\subseteq$ K. Suppose that A/K is nilpotent and let B/K be the unuque minimal ideal of A/K. Now A/K=Nil(A/K)=B/K by Theorem 2.4 of [4] and dim(A/K)=1, a contradiction. Thus A/K is not nilpotent.

\vspace{5mm}

Proposition 19 Let A be a solvable Leibniz algebra and let K be a primitive ideal in A. Then K is generalized Frattini in A if and only if K is a proper subalgebra of Nil(A).

\vspace{5mm}

Proof. Suppose that K is generalized Frattini in A and let B/K be the unique minimal ideal of A/K. Since A is solvable, B/K=Nil(A/K) using Theorem 2.4  of [4]. Then B/K=Nil(A/K)=NilA)/K by Proposition 13. Hence Nil(A)=B and K is a proper subalgebra of Nil(A).

\vspace{5mm}

Conversely, let the ideal K be a proper subalgebra of Nil(A). Then Nil(A/K)=B/K. Let H be an ideal in A such that K is properly contained in H and H/K is nilpotent. Then H/K $\subseteq$ Nil(A/K)=Nil(A)/K. Hence H $\subseteq$ Nil(A) and H is nilpotent. Therefore K is generalized Frattini in A by Theorem 7.

\vspace{5mm}

Theorem 20. Let A be a solvable Leibniz algebra and let K be a primitive ideal in A. Let B/K be the unique minimal ideal in A/K. Then K is generalized Frattini in A if and only if B=Nil(A).

\vspace{5mm}

Proof. Suppose that K is generalized Frattini in A. Then Nil(A/K)=Nil(A)/K by Theorem 7 . Since A is solvable, B/K=Nil(A/K)=Nil(A)/K and B=Nil(A).

\vspace{5mm}

Conversely, suppose that B=Nil(A). Since A ios solvable, Nil(A/K)=B/K=Nil(A)/K and K is a proper subalgebra of Nil(A). If  N is any ideal of A with N/K nilpotent, then N is nilpotent. Hence K is generalized Frattini in A by Theorem 7.

\vspace{5mm}

Corollary 21. Let K be a primitive ideal of a solvable Leibniz algebra  A. If  K is generalized Frattini in A, then K is maximal with respect to the generalized Frattini property in A.

\vspace{5mm}

Proof. Suppose that H is generalized Frattini ideal in A such that K $\subseteq$ H. Then H is nilpotent by Proposition 1 and hence, H is contained in Nil(A). Let B/K be the unique minimal ideal in A/K. By Theorem 20, B=Nil(A). Hence, either H=K or H=Nil(A) since K $\subseteq$ H $\subseteq$ Nil(A).Suppose that H=Nil(A). Then, by Proposition 10, every solvable ideal of A is nilpotent. In particular, A is nilpotent.This contradicts Lemma 18. Thus H=K and K is maximal with respect to the generalized Frattini property.

\vspace{5mm}

IV. INTERSECTIONS OF CERTAIN MAXIMAL SUBALGEBRAS

\vspace{5mm}

Let R(A) be the intersection of all maximal subalgebras that are ideals in A and T(A) be the intersection of all maximal subalgebras that are not ideals in A. As for F(A), T(A) may not be an ideal and we let $\tau$(A) be the largest ideal of A contained in T(A). Of course $\Phi$(A)=R(A)$\cap \tau$(A).

\vspace{5mm}

An algebra is power solvable if all subalgebras generated by a single element are solvable. For Leibniz algebras, these subalgebras are cyclic subalgebras in which A$^2$=Leib(A) is abelian (Section 4 of [5]). Hence Leibniz algebras are power solvable and the following result, Theorem 2.8 of [15] holds.

\vspace{5mm}

Lemma 22. If $\Phi$(A)=0, then $\tau$(A)= Z(A)=Z$^*$(A) where Z$^*$(A) is the final term in the upper central series of A.

\vspace{5mm}

Proposition 23.  $\tau$(A) is generalized Frattini in A.

\vspace{5mm}

Proof. By the last lemma, $\tau$(A)/$\Phi$(A)=Z(A)/$\Phi$(A)=Z(A/$\Phi$(A)). By Theorem 1, Z(A/$\Phi$(A)) is generalized Frattini in A/$\Phi$(A). Thus $\tau$(A) is generalized Frattini in A

\vspace{5mm}

Proposition 24. Let A be a non-nilpotent Leibniz algebra with $\Phi$(A)=0. Then any ideal H that is a maximal generalized Frattini ideal in A  contains $\tau$(A).

\vspace{5mm}

Proof.  H+$\Phi$(A) is generalized Frattini by Theorem 1. Hence H+$\Phi$(A)=H since H is maximal generalized Frattini. Thus $\Phi$(A) $\subseteq$ H. Then (H+$\tau$(A))/$\Phi$(A)=H/$\Phi$(A)+Z(A/$\Phi(A)$) which is generalized Frattini in A/$\Phi$(A) by Theorem 1.  . Hence H+Z(A)=H+$\tau$(A) which is generalized Frattini by Theorem 1  . Using maximality of H, H=H+$\tau$(A). Hence $\tau$(A) $\subseteq$H.

\vspace{5mm}

Proposition 25. Let A be a non-nilpotent Leibniz algebra with $\Phi$(A)=0. Then any ideal H that is maximal with respect to the generalized Frattini property in A contains Z$^*$(A).

\vspace{5mm}

Proof. (Z$^*$(A)+$\Phi$(A))/$\Phi$(A) is contained in Z$^*$(A/$\Phi$(A))=$\tau$(A)/$\Phi$(A). Hence Z$^*$(A) $\subseteq \tau$(A) $\subseteq$H by Proposition 24.

\vspace{5mm}

Theorem 26. A is nilpotent if and only if R(A) $\subseteq \tau$(A). [2] or [10]

\vspace{5mm}

Proof. If A is nilpotent, then all maximal subalgebras are ideals in A [10]. Hence $\tau$(A)=A and R(A) $\subseteq \tau$(A). Conversely, if R(A) $\subseteq \tau$(A), then $\Phi$(A)=R(A) which contains A$^2$ by Lemma 2.3 of [15]. Hence all maximal subalgebras are ideals and A is nilpotent.

\vspace{5mm}

Corollary 27. A is nilpotent if and only if $\Phi$(A)=A$^2$.

\vspace{5mm}

Proof. If A is nilpotent, then all maximal subalgebras are ideals. Hence $\tau$(A)=A and the result follows from Theorem 26. Conversely, if the condition holds, then all maximal subalgbras are ideals and A is nilpotent by [2].

\vspace{5mm}

Lemma 28. Let N be an ideal in A.

1. F(A)+N/N $\subseteq$ F(A/N)

2. If N $\subseteq$ F(A), then F(A)/N=F(A/N)

\vspace{5mm}

Proof. 1. If M/N is a maximal subalgebra of A/N, then M is a maximal subalgebra of A. Thus F(A) $\subseteq \cap_{M/N maximal in A/N}$M. Hence F(A)+N/N $\subseteq$ F(A/N).

\vspace{5mm}

2. N is contained in all maximal subalgebras of A. M/N is a maximal subalgebra of A/N if and only if M is a maximal subalgebra of A. Then F(A)/N=F(A/N).

\vspace{5mm}

Let nFrat(A) be the intersection of all maximal ideals of A. Then following the same arguements as in the last lemma, we obtain

\vspace{5mm}

Lemma 29. Let N be an ideal in A.

1. nFrat(A)+N/N $\subseteq$ nFrat (A/N)

2. If N $\subseteq$ nFrat(A), then nFrat(A)/N=nFrat(A/N).

\vspace{5mm}

In Lie algebras, the Frattini ideal is nilpotent. However, nFrat(A) and R(A) need not be nilpotent. The same results hold for Leibniz algebras also. For x in A, L$_x$ denotes left multiplication by x and A$_0$(x) and A$_1$(x) are the corresponding Fitting components for L$_x$. A$_0$(x) is the Engel subalgebra from section II.

\vspace{5mm}

Proposition 30. $\Phi$(A) is nilpotent,  but nFrat(A) and R(A) need not be nilpotent.

\vspace{5mm}

Proof. Let x$\in \Phi$(A). Then A$_1$(x) $\subseteq \Phi$(A) since $\Phi$(A) is an ideal in A. Hence A$_0$(x), a subalgebra of A, supplements $\Phi$(A). Hence A$_0$(x)=A and L$_x$ is nilpotent for all x in $\Phi${A) and $\Phi$(A) is nilpotent by Engel's theorem. For the examples, let A=gl(p,F) where F has characteristic p. The Z(A) $\subset$ A$^2$=sl(p,F) which is the only maximal ideal of A and nFrat(A) is not nilpotent. Since sl(p,F) is the only maximal subalgebra that is an ideal, R(A)=sl(p,F) is not nilpotent.

\vspace{5mm}

Both $\Phi$(A) and nFrat(A) are contained in R(A). We find other results of this type.

\vspace{5mm}

Proposition 31. $\Phi$(A) $\subseteq$ nFrat(A) $\subseteq$ R(A).
\vspace{5mm}

Proof If N is a maximal ideal of A, then $\Phi$(A)+N can not equal A, hence $\Phi$(A) $\subseteq$ N. Hence $\Phi$(A) $\subseteq$ nFrat(A). The other inclusion is clear.

\vspace{5mm}

Proposition 32. If A is solvable, then R(A)=nFrat(A).

\vspace{5mm}

Proof. If N is a maximal ideal of a solvable Leibniz algebra, then dim (A/N)=1. Thus every maximal ideal is a maximal subalgebra that is an ideal. The converse also holds. Hence the result holds.

\vspace{5mm}

If A is not solvable, then Proposition 32 does not hold.

\vspace{5mm}

Example 33. Let A=sl(2,F). Since A is simple, nFrat(A)=0. Since A contains no maximal subalgebras that are ideals, A=R(A).

\vspace{5mm}

Theorem 34. A is nilpotent if and only if $\Phi$(A)=nFrat(A)=R(A).

\vspace{5mm}

Proof. If A is nilpotent, then all maximal subalgebras are ideals and R(A)=$\Phi$(A) and Proposition 31 gives the result.

\vspace{5mm}

Suppose that the three ideals are equal. Since $\tau$(A) always contains $\Phi$(A), R(A) $\subseteq \tau$(A). Then A is nilpotent by Theorem 26.

\vspace{5mm}

V. NON GENERATORS

\vspace{5mm}

Following results in groups and Lie algebras, we give Leibniz algebra characterizations of F(A), R(A) and nFrat(A) in terms of non-generators. 

\vspace{5mm}

A subset S of A that is closed under multiplications by elements of A is called a normal set in A. An element x $\in$ A is a normal non-generator in A if whenever A= $<x,T>$ it follows that A=$<T>$.

\vspace{5mm}

Proposition 35. F(A) consists of the non-generators of A.

\vspace{5mm}

Proof. Let x$\in$ F(A) and let A=$<x,H>$ where H is a subalgebra of A. If H $\neq$ A, then H $\subseteq$ M for some maximal subalgebra of A. Then $<H,x> \subseteq$ M, a contradiction. Hence A=H and x is a non-generator of A. Suppose that x is not in F(A). Let M be a maximal subalgebra of A which does not contain x. Then M $\subset <x,M>  $=A. But M $\neq$ A, so x is a not a non-generator of A.

\vspace{5mm}

Proposition 36. R(A) equals the set of normal non generators of A.

\vspace{5mm}

Proof. Suppose that x $in$ R(A) and A=$<x,S>$ for a normal subset of S of A. Is $<S> \neq$ A, then dim A= dim $<S>$+1, so $<S>$ is a maximal subalgebra of A that is an ideal in A. This contradicts that x$\in$ R(A), so $<S>$=A, and x is a normal non-generator of A.

\vspace{5mm}

Conversely, suppose that x is not in R(A). Then x is not in a maximal subalgebra, M, that is an ideal in A. M is a normal subset of A and A=$<x,M>$ but A $\neq$ M. Hence x is not a noral non-generator in A.

\vspace{5mm}

Let X be a subset of A. Let X$^A$ be the smallest ideal in A that contains X. An element x $\in$ A is called an n-nongenerator if for every subset X of A, A=X$^A$ whenever A=$<x,X>^A$

\vspace{5mm}

Lemma 37. For any x$\in$ A and any subset X of A, $<x,X>^A$=$<x^A,X^A>$=$x^A$+X$^A$.

\vspace{5mm}

Proof. Both x$^A$ and X$^A$ are contained in $<x,X>^A$. Therefore $<x^A,X^A> \subseteq <x,X>^A$ and $x^A+X^A \subseteq <x,X>^A$. Since $<x,X> \subseteq <x^A,X^A>$, it follows that $<x,X>^A \subseteq <x^A,X^A>$. Also $x^A, X^A \subseteq x^A+X^A$. Thus $<x^A,X^A> \subseteq x^A+X^A$.

\vspace{5mm}

Proposition 38. nFrat(A) is the set of n-nongenerators for A.

\vspace{5mm}

Proof. let T be the collection of all n-nongenerators of A . Suppose x is in T but not in nFrat(A). Let N be a maximal ideal such that x is not in N. Then x$^A$+N=A. Hence $<x,N>^A$=A. Thus N=N$^A$=A, a contradiction. Hence T $\subseteq$ nFrat(A).

\vspace{5mm}

Conversely, let x be in nFrat(A) but not in T. There exists a subset S of A such that A=$<x,S>^A$ but S$^A$ is properly contained in A. Therefore S$^A$ is a proper ideal of A and x is not in S$^A$. By Lemma 37, A=$<x,S>^A=x^A+S^A$. Let M be maximal with respect to the properties for S$^A$: x is not in M, M is an ideal of A, S$^A \subseteq$ M and A=x$^A$+M. We claim that M is a maximal ideal of A. If not, let N be an ideal properly between M and A. Then A=x$^A$+M=x$^A$+N. By the maximality conditions on M, x$\in$ N. Therefore, x$^A \subseteq$ N and A=N, a contradiction. Thus M is a maximal ideal in A. Since x is not in M, it follows that x is not in nFrat(A), a contradiction. Hence, whenever A=$<x,S>^A$, it follows that A=$S^A$ and x is an n-nongenerator of A. Hence nFrat(A) $\subseteq$ T and the result holds.

\vspace{5mm}

References

\vspace{5mm}

Sh. Ayupov, B. Omirov, On Leibniz algebras, In Algebra and Operator Theory, Colloquium in Tashkent, Kluwer, 1998

D. Barnes, Some theorems on Leibniz algebras, Comm in Alg 39, 2463-2472 (2011)

C. Batten, L. Bosko-Dunbar, A. Hedges, J.T. Hird, K. Stagg, E. Stitzinger, A Frattini Theory for Leibniz algebras, Comm, in Algebra 41,1547-1557 (2013)

C Batten Ray, A. Combs, N Gin, A. Hedges, J. Hird, L. Zack, Nilpotent Lie and Leibniz algebras, Comm. in Algebra 42, 2404-2410 (2014)

C Batten Ray, A. Hedges, E. Stitzinger, , Classifying several classes of Leibniz algebras, Alg. and Rep. Theory 17, 703-712 (2014)

J. Beidleman and T. Seo, Generalized Frattini subgroups of finite groups, Pacific J. Math.,  23, 441-450, (1967)

I. Demir, K Misra, E Stitzinger, On some stuctures of Leibniz algebras, A.M.S. Cotemporary Math. 623,  41-54 (2014)

J-L Loday, Une version non commutative des algebres de Lie: les algebres de Leibniz, Enseign Math. 39, 269-293 (1993)

K Stagg, Analogues of the Frattini subalgebra, International Electronic Journal of Algebra 9, 124-132 (2011)

D. Towers, A Frattini theory for algebras, Proc. Lond. Math. Soc, 27, 440-462 (1973)

D. Towers, Two ideals of an algebra closely related to it's Frattini ideal. Archiv Math. 35 112-120 (1980)

\end{document}